\documentclass[11pt]{article}
\usepackage{amssymb}
\usepackage{graphicx}
\usepackage{amsmath}
\usepackage{amsfonts}
\usepackage[english]{babel}
\usepackage{amsthm}
\usepackage{latexsym}
\usepackage{hyperref}
\newtheorem{theorem}{Theorem}
\newtheorem{corollary}{Corollary}

\newtheorem{problem}{Problem}
\newtheorem{remark}{Remark}
\def\neweq#1{\begin{equation}\label{#1}}
\def\endeq{\end{equation}}
\def\eq#1{(\ref{#1})}

\newcommand{\cC}{{\mathcal C}}
\newcommand{\M}{{\mathcal M}}

\newcommand{\Sc}{{\mathcal S}}
\newcommand{\N}{\mathbb{N}}
\newcommand{\LL}{{\mathcal L}\,}

\newcommand{\R}{\mathbb{R}}

\newcommand{\integ}{\int_{\R^n}}

\newcommand{\rhoa}{{ a }}

\newcommand{\dist}{{\rm dist}}

\newcommand{\C}{{\mathbb{C}}}

\newcommand{\cE}{{\cal E}}

\newcommand{\ab}{\alpha,\beta}
\newcommand{\subab}{\substack{|\alpha|=m \\ |\beta|=m}}
\newcommand{\subabl}{\substack{|\alpha|\leq m \\ |\beta|\leq m}}

\newcommand{\darr}[4]{{\left\{\begin{array}{ll}
   {#1}&{#2}\\
   {#3}&{#4}
 \end{array}\right.}}
 
\newcommand{\darrn}[4]{{\begin{array}{ll}
   {#1}&{#2}\\[0.2cm]
   {#3}&{#4}
\end{array}}}

\newcommand{\ia}{({\rm i})}
\newcommand{\ib}{({\rm ii})}

\newcommand{\dom}{{\rm Dom}}
\newcommand{\comb}[2]{c^{#1}_{#2}}

\begin{document}
\title{Higher order linear parabolic equations}

\author{Gerassimos BARBATIS\\
{\small Department of Mathematics}\\
{\small University of Athens}\\
{\small Panepistimioupolis}\\
{\small 15784 Athens (Greece)}
\and Filippo GAZZOLA \\
{\small Dipartimento di Matematica}\\
{\small Politecnico di Milano}\\
{\small Piazza Leonardo da Vinci 32}\\
{\small 20133 Milano (Italy)}}

\date{{\em Dedicated to Patrizia Pucci, on the occasion of her 60th birthday.}}

\maketitle

\begin{abstract}
We first highlight the main differences between second order and higher order linear parabolic equations. Then
we survey existing results for the latter, in particular by analyzing the behavior of the convolution kernels.
We illustrate the updated state of art and we suggest several open problems.\par
\end{abstract}

\section{Introduction}

The Cauchy problem in $\R^n$ ($n\ge1$) for higher order ($m\ge2$) linear parabolic equations
\neweq{general}
\left\{
\begin{array}{lll}
\displaystyle u_t+(-1)^m\sum_{\subabl} D^\alpha\{a_{\alpha,\beta}D^\beta\}u=0\qquad & \mbox{in }\R^n\times\R_+ \, ,\\
\\
u\left( x,0\right) =u_{0}(x)\qquad & \mbox{in }\mathbb{R}^{n} \, ,
\end{array}
\right.
\endeq
has recently attracted some interest, due to its somehow surprising and unexpected properties, strikingly different when compared with the corresponding second order parabolic
equation, that is, when $m=1$. The purpose of the present paper is to survey existing results about problem \eq{general}
and to suggest several open problems whose solution would contribute towards the formation of a complete theory.\par
Even in the simplest situation when \eq{general} becomes the polyharmonic heat equation
\begin{equation}\label{linear2}
\left\{
\begin{array}{ll}
u_{t}+(-\Delta)^{m}u=0\qquad & \mbox{in }\R^n\times\R_+ \, , \\
u\left( x,0\right) =u_{0}(x)\qquad & \mbox{in }\mathbb{R}^{n}\ ,
\end{array}
\right.
\end{equation}
important differences appear and many questions are still open. As was first observed by Evgrafov-Postnikov \cite{ep}, the kernels of the heat operators
in \eq{linear2} depend on the space dimension, contrary to the classical second order heat operator; this apparently harmless fact, already claims a lot
of work in order to obtain fine
qualitative properties of the solution to \eq{linear2}. When $u_{0}\in C^{0}\cap L^{\infty }\left( \mathbb{R}^{n}\right)$,
problem (\ref{linear2}) admits a unique global in time bounded solution explicitly given by
\begin{equation}\label{fmn}
u(x,t)=\alpha t^{-n/2m}\int_{\mathbb{R}^{n}}u_{0}(x-y)f_{m,n}\Big( \frac{|y|}{t^{1/2m}}\Big) \,dy\, , \qquad (x,t)\in \R^n\times\R_+\ ,
\end{equation}
where $\alpha=\alpha_{m,n}>0$ is a suitable normalization constant and
\neweq{polyker}
f_{m,n}(\eta)=\eta^{1-n}\int_0^\infty e^{-s^{2m}}(\eta s)^{n/2}J_{(n-2)/2}(\eta s)\, ds\ ,
\endeq
see \cite{carinhas}. Here and below, $J_{\nu}$ denotes the $\nu$-th Bessel function. So, not only the kernels $f_{m,n}$ depend on $n$, but also they are not available in a simple
form. Due to the presence of Bessels functions in \eq{fmn}, the solution to \eq{linear2} exhibits oscillations and this fact has two main consequences. First, the positivity
preserving property fails; it is in general {\em false} that positivity of the initial datum $u_0$ yields positivity of the solution $u$.
Second, in order to prove global existence or finite time blow-up for corresponding semilinear equations, comparison principles cannot be used; for this reason,
Galaktionov-Poho\v{z}aev \cite{galaktionovpohozaev} introduced a new method based on {\em majorizing order-preserving operators} which, basically, consists in taking
the convolution of the initial datum $u_0$ with the absolute value of the kernel $f_{m,n}$.\par
The asymptotic behavior of the solution to the second order heat equation can be described with some precision also thanks to the so-called Fokker-Plank equation
obtained by exploiting the self-similar structure of the fundamental solution. But the Fokker-Plank operator corresponding to \eq{linear2} is not self-adjoint if $m\ge2$
and this brings several difficulties to the analysis of its spectral properties; these difficulties were partially overcome in a fundamental paper by
Egorov-Galaktionov-Kondratiev-Poho\v zaev \cite{egkp2}. However, most of the classical methods usually exploited for the second order heat equation do not apply.
For instance, any reasonable Lyapunov functional becomes very complicated due to the presence of higher order derivatives, too many terms appear and the
study of their signs is out of reach. Also standard entropy methods fail, due to the change of sign of the kernels $f_{m,n}$: the second order entropy is $\int u\log u$
and cannot be considered because the solution $u$ to \eq{linear2} changes sign also for positive data. The sign change of the kernels also forbids to analyze the
behavior of suitable scaled ratios such as $u/f_{m,n}$ in order to obtain Ornstein-Uhlenbeck-type equations.\par

The fact that the functions $f_{m,n}$ exhibit oscillations also implies that the semigroup associated to \eq{general} is not Markovian if $m\ge2$; this yields important complications in extending the $L^2$ theory to an $L^p$ theory. In the second order case one uses the Markovian properties of the $L^2$ semigroup to prove that it extends to a contraction semigroup in $L^p$. This then leads to heat kernel estimates, a topic extensively studied in the past 25 years.
For $m\geq 2$ and $L^{\infty}$ coefficients
the situation is reversed: one first obtains heat kernel estimates and then applies them in order to develop the $L^p$ theory. The heat kernel estimates depend essentially upon the validity of the Sobolev embedding $H^m(\R^n) \subset C^0(\R^n)$, hence an important distinction arises depending on the dimension $n$. This is in contrast to the second order case where the theory does not depend on such an embedding.
\par
The problem of obtaining sharp heat kernel estimates is itself very interesting. To put it into context, one needs to go back to short time asymptotic estimates, first proved by Evgrafov-Postnikov \cite{ep} for constant coefficient equations and later extended by Tintarev \cite{tintarev} for variable smooth coefficients. Progress has been made in the past years in obtaining sharp heat kernel bounds, but several important questions remain open.
\par
Further recent results are available for \eq{linear2}. In \cite{berchio} the positivity preserving property
is studied in presence of a source $f(x,t)$.
In \cite{chol1} the solvability of the Cauchy problem \eq{linear2} (with $m=2$) in presence of an irregular datum $u_0$ is studied and the presence of a strongly continuous
analytic semigroup is proved. Finally, we mention that more general linear problems were considered in  \cite{quesada} whereas
the stability method for higher order equations was studied in \cite[Chapter 12]{galaktionovvazquez}.\par
For the above reasons, many natural questions arise. In this paper, we mainly focus our interest on the fundamental solution (heat kernel) of (\ref{general}) and on
positivity preserving property (ppp from now on). As already mentioned, these problems are by now very well understood in the second order case where heat kernels
have been extensively studied in very general frameworks, while ppp holds as consequence of the positivity of the Gaussian heat kernel (maximum principle).
In the higher order case the situation is considerably more complicated and it is precisely our purpose to give an updated state of art as well as a number of
open problems still to be solved in order to reach a satisfactory theory.\par
In Section \ref{section:hke} we study various properties of the heat kernel of the general problem \eq{general}. We avoid any local regularity assumptions on the coefficients, and we start with Davies' results \cite{davies1,davies2,davies3,davies4} on operators with $L^{\infty}$ coefficients, omitting reference
to earlier work where local regularity assumptions were imposed. We then tackle
the $L^p$ theory, emphasizing the dimensional dependence.
We proceed to present the short time asymptotics of Evgrafov-Postnikov and Tintarev and sharp heat kernel estimates, including results on non-uniformly elliptic operators.
In the last part of Section \ref{section:hke} we restrict our attention to constant coefficients case, namely equation \eq{linear2}. In this simplified
situation, especially if $m=2$, much more can be said on the behavior of the kernels; in particular, we exhibit fine properties of their moments.\par
In Section \ref{fokker} we transform \eq{linear2} into a Fokker-Planck-type equation and we recall an important result by
Egorov-Galaktionov-Kondratiev-Poho\v zaev \cite{egkp2} about the spectrum of the corresponding (non self-adjoint) operator. In Section \ref{asymptoticbeh} we determine the
behavior of the moments of the solution to the Fokker-Planck equation in the fourth order case $m=2$.\par
In Section \ref{ppp} we recall the results which describe the way how the ppp may fail and we discuss the possibility of finding a limit decay of the datum $u_0$ for which
ppp may still hold.

\section{Heat kernel estimates}\label{section:hke}

In this section we survey some properties of the heat kernel of problem (\ref{general}).
We first discuss the case where the operator has $L^{\infty}$
coefficients, then we extend some results to the ``singular case'' where the coefficients are merely assumed to be in $L^{\infty}_{\rm loc}$, finally
we specialize to the the simplest case of constant coefficients for fourth order equations: the biharmonic heat kernel.

\subsection{Semigroup generation}

Problem (\ref{general}) is to be understood in the $L^2$-sense, and for this we need to properly define the elliptic operator
\[
(Hu)(x)=(-1)^m\sum_{\subabl}D^{\alpha}\{a_{\ab}(x)D^{\beta}u\}
\]
as a self-adjoint operator in $L^2(\R^n)$. For this we start with real-valued functions
$a_{\ab}(x)=a_{\beta,\alpha}(x)$, $|\alpha|,|\beta|\leq m$, in $L^{\infty}(\R^n)$ and we define the quadratic form
\[
Q(u)=\int_{\R^n}\sum_{\subabl}a_{\ab}(x)D^{\alpha}u D^{\beta}\bar{u}\, dx
\]
on $\dom(Q)= H^m(\R^n)$. Our main ellipticity assumption is that G{\aa}rding's inequality
\begin{equation}
Q(u)\geq c_1\|u\|_{H^m(\R^n)}^2 -c_2 \|u\|_{L^2(\R^n)}^2  \; , \quad u\in H^m(\R^n),
\label{barb:garding}
\end{equation}
is satisfied for some $c_1,c_2>0$. It then follows that the form $Q$ is closed; the operator $H$ is defined as the self-adjoint operator on $L^2(\R^n)$
associated to the quadratic form $Q$.
It is well-known \cite[Theorem 7.12]{agmon} that inequality (\ref{barb:garding}) implies that the principal symbol of $H$ satisfies
\[
\sum_{\subab}a_{\ab}(x)\xi^{\alpha+\beta} \geq c_1|\xi|^{2m}, \; \quad \xi\in\R^n , \; x\in\R^n
\]
and that the converse implication is true for uniformly continuous coefficients.

We first consider the question of existence of a heat kernel together with pointwise estimates. The heat kernel $K(t,x,y)$ of $H$ is, by definition, the integral kernel of the semigroup $e^{-Ht}$, provided such kernel exists.
Hence it represents the solution $u(x,t)$ of \eq{general} in the sense that
\[
u(x,t)=\int_{\R^n}K(t,x,y)u_0(y)dy \, , \quad \forall (x,t)\in\R^n\times\R_+ \, .
\]
The results depend on whether the order $2m$ of $H$ exceeds or not the dimension $n$.
\begin{theorem}\label{thm:barb:1}{\rm \cite[Lemma 19]{davies1}, \cite[Theorem 1.1]{ElRo}, \cite[Proposition 28]{AuTc}}\par
If $2m\geq n$ then the semigroup $e^{-Ht}$ has a continuous integral kernel $K(t,x,y)$. Moreover there exist positive constants $c_i$, $i=1,2,3$, such that
\begin{equation}
\label{barb:ka:10}
| K(t,x,y)| < c_1t^{-\frac{n}{2m}}\exp\left\{ -c_2\frac{|x-y|^{\frac{2m}{2m-1}}}{t^{\frac{1}{2m-1}}} + c_3t   \right\},
\end{equation}
for all $t\in\R_+$ and $x,y\in\R^n$.
\end{theorem}

One application of this theorem is the extension of the $L^2$-theory to $L^p(\R^n)$. The fact that the semigroup $e^{-Ht}$ is not Markovian makes
this problem quite different from the second order case $m=1$.

\begin{theorem}\label{thm:barb:3}{\rm \cite[Theorems 20 and 21]{davies1}}\par
Assume that $2m\geq n$. The semigroup $e^{-Hz}$, ${\rm Re}\, z>0$, extends from $L^2(\R^n)\cap L^p(\R^n)$ to a bounded holomorphic semigroup $T_p(z)$ on $L^p(\R^n)$ for all $1\leq p\leq\infty$. Moreover, for $1\leq p<\infty$ the semigroup $T_p(z)$ is strongly continuous and its generator $-H_p$ has spectrum which is independent of $p$.
\end{theorem}

In the case $2m<n$ critical Sobolev embedding into $L^p$ spaces appear and the situation is different.

\begin{theorem}\label{thm:barb:4}{\rm \cite[Theorem 10]{davies2}}\par
Assume that $2m<n$. Let $p_c=2n/(n-2m)$ be the Sobolev exponent and let $q_c=2n/(n+2m)$ denote its conjugate.\par
$\ia$ The semigroup $e^{-Hz}$ extends to a strongly continuous bounded holomorphic semigroup $T_p(z)$ on $L^p(\R^n)$ for
all $q_c\leq p\leq p_c$. Moreover the spectrum of the generator $-H_p$ of $T_p(z)$ is independent of $p$.\par
$\ib$ Assume that $m$ is even. For $p\not\in [q_c , p_c]$ there exists an operator $H$ of the above type for which the operator $e^{-Ht}$ does not extend from $L^2(\R^n)\cap L^p(\R^n)$ to a bounded operator on $L^p(\R^n)$, for any $t>0$. In particular the semigroup $e^{-Ht}$ does not have an integral kernel satisfying (\ref{barb:ka:10}).
\end{theorem}
We note that when $m$ is odd a result analogous to (ii) is valid for elliptic systems \cite{davies2}. We also note that if the coefficients are sufficiently regular
then a Gaussian heat kernel estimate is valid without any restriction on the dimension; see \cite{davies3} and references therein for more details.

\subsection{Short time asymptotic estimates}

In this subsection we make the additional assumption that the coefficients $\{a_{\ab}(x)\}$ are smooth.
We consider the problem \eq{general} and denote by
\[
A(x,\xi)= \sum_{\subab}a_{\ab}(x)\xi^{\alpha+\beta} \; ,
\]
the corresponding principal symbol, which satisfies
\[
 c^{-1}|\xi|^{2m}\leq A(x,\xi) \leq c|\xi|^{2m} \; , \qquad \xi\in\R^n , \; x\in\R^n,\\[0.2cm]
\]
for some $c>0$.
The following notion of {\em strong convexity} was first introduced by Evgrafov-Postnikov \cite{ep}. For a multi-index $\gamma$ with $|\gamma|=2m$
we denote $\comb{2m}{\gamma} =(2m)! /(\gamma_1 ! \ldots \gamma_n !)$.
We define the functions $b_{\gamma}(x)$, $|\gamma|=2m$, by requiring that
\[
A(x,\xi) =\sum_{|\gamma|=2m}\comb{2m}{\gamma}b_{\gamma}(x)\xi^{\gamma} \; , \quad \xi\in\R^n , \; x\in\R^n.
\]
{\bf Definition.} The symbol $A(x,\xi)$ is strongly convex if the quadratic form
\[
\Gamma(x,v)=\sum_{\subab}b_{\alpha+\beta}(x)v_{\alpha}\bar{v}_{\beta} \; , \quad v=(v_{\alpha})\in\C^{\nu},
\]
is positive semi-definite for all $x\in\R^n$. It is known \cite[Section 1]{ep} that strong convexity implies that the matrix
$\{A_{\xi_i\xi_j}(x,\xi)\}_{i,j}$ is positive definite for all $x\in\R^n$ and $\xi\in\R^n\setminus\{0\}$.

We first consider operators with constant coefficients so that $K(t,x,y)=K(t,x-y,0)$. We set
\begin{equation}
\label{sigma}
\sigma_m =(2m-1)(2m)^{-\frac{2m}{2m-1}}\sin\Big(\frac{\pi}{4m-2}\Big).
\end{equation}

\begin{theorem}\label{barb:thm:40}{\rm \cite[Theorem 4.1]{ep}}\par
Assume that $H$ is homogeneous of order $2m$ with constant coefficients and that the symbol $A(\xi)$ is strongly convex. Let
\begin{equation}\label{eq:barb:200}
p(\xi)=\max_{{\substack{\eta\in\R^n \\ \eta\neq 0}}}\frac{\xi\cdot\eta}{A(\eta)^{1/2m}}, \qquad \xi\in\R^n.
\end{equation}
There exists a positive function $S(x)$ such that for any $x\in\R^n$, $x\neq 0$, we have
\begin{eqnarray*}
K(t,x,0)&=& S(x)t^{-\frac{n}{2(2m-1)}}\cos\Big( \sigma_m\frac{ p(x)^{\frac{2m}{2m-1}}}{t^{\frac{1}{2m-1}}} \cot\Big(\frac{\pi}{4m-2}\Big) -\frac{n(m-1)}{4m-2} +o(1)
\Big)\\
&& \qquad \times \exp\left\{ -\sigma_m\frac{ p(x)^{\frac{2m}{2m-1}}}{t^{\frac{1}{2m-1}}}  \right\}
\end{eqnarray*}
as $t\to 0$.
\end{theorem}
In order not to become too technical we refer to \cite{ep} for the precise definition of $S(x)$; we note however that it is positively homogeneous of degree $-n(m-1)/(2m-1)$ in $x\in\R^n\setminus\{0\}$.

To extend Theorem \ref{barb:thm:40} to the case of variable smooth coefficients we need some elementary notions of Finsler geometry.
Very roughly, one can say that a Finsler metric is the assignment of a norm at each tangent space of a manifold.
In our context, extending (\ref{eq:barb:200}) we define
\[
p(x,\xi)=\max_{{\substack{\eta\in\R^n \\ \eta\neq 0}}}\frac{\xi\cdot\eta}{A(x,\eta)^{1/2m}}\; , \qquad x\in\R^n , \; \xi\in\R^n.
\]
This defines a Finsler metric on $\R^n$ in the sense that
\neweq{defFin}
\mbox{$p(x,\xi)=0$ if and only if $\xi=0$\qquad and\qquad$p(x,\lambda\xi)=|\lambda|p(x,\xi)$, $\lambda\in\R$.}
\endeq
In Finsler geometry the definition is typically complemented by
\neweq{Fin2}
\mbox{the matrix $\{g_{ij}\}:= \frac{1}{2}\frac{\partial^2p(x,\xi)^2}{\partial\xi_i\partial\xi_j}$ is positive definite for all
$x\in\R^n$ and $\xi\in\R^n\setminus\{0\}$.}
\endeq
For our purposes we shall not assume \eq{Fin2} except in Theorem \ref{barb:thm:approx_fm} below.
We note however that if \eq{defFin}-\eq{Fin2} are valid then the map $\xi\mapsto p(x,\xi)$ is indeed a norm for all $x\in\R^n$.
The length of an absolutely continuous path, $\gamma=\gamma(t)$, $0\leq t\leq 1$, is then defined as
\begin{equation}\label{barb:def:fd}
l(\gamma)=\int_0^1 p(\gamma(t) ,\dot{\gamma}(t))dt \; ,
\end{equation}
and the Finsler distance between two points $x,y\in\R^n$ is given by
\[
d(x,y)=\inf \{l (\gamma) \; : \; \mbox{$\gamma$ has endpoints $x$ and $y$}\}.
\]
\begin{theorem}\label{barb:thm:30}{\rm \cite[Theorem 1.1]{tintarev}}\par
Assume that the operator $H$ is homogeneous of order $2m$ with smooth coefficients and that the principal symbol $A(x,\xi)$ is strongly convex. Assume
further that the matrices $\{a_{\ab}(x)\}_{|\alpha|=|\beta|=m}$ and $\{A_{\xi_i\xi_j}(x,\xi)\}_{1\leq i,j\leq n}$ are both positive definite uniformly
in $x\in\R^n$ and $\xi\in S^{n-1}$. Then there exist functions $v_k(t,x,y)$, $k=0,1,\ldots$, such that the following is true:
for any $x\in\R^n$ there exists $\delta>0$ such that for $0<|x-y|<\delta$ the following asymptotic expansion is valid as $t\to 0$:
\begin{equation}
K(t,x,y) \sim \sum_{k=0}^{\infty}t^{\frac{k-\frac{n}{2}}{2m-1}}v_k(t,x,y)\exp\left\{ -\sigma_m \frac{d(x,y)^{\frac{2m}{2m-1}}}{t^{\frac{1}{2m-1}}} \right\}.
\label{barb:5}
\end{equation}
The functions $v_k(t,x,y)$ oscillate and are bounded and smooth with respect to $t$.
\end{theorem}
Estimate (\ref{barb:5}) is meant in the sense that for each $N\geq 1$ and for small enough $t>0$ there holds
\[
\bigg| K(t,x,y) - \sum_{k=0}^{N}t^{\frac{j-\frac{n}{2}}{2m-1}}v_k(t,x,y)
\exp\left\{ -\sigma_m \frac{d(x,y)^{\frac{2m}{2m-1}}}{t^{\frac{1}{2m-1}}} \right\}\bigg|
\leq c_N t^{\frac{N+1-\frac{n}{2}}{2m-1}}\exp\left\{ -\sigma_m \frac{d(x,y)^{\frac{2m}{2m-1}}}{t^{\frac{1}{2m-1}}} \right\}.
\]

\subsection{Sharp heat kernel bounds}

We now return to the general framework of operators with $L^{\infty}$ coefficients satisfying G{\aa}rding's inequality (\ref{barb:garding}).
We assume that $2m>n$ so that the heat kernel estimate (\ref{barb:ka:10}) is valid and we present certain theorems that provide additional information
on the constant $c_2$ in (\ref{barb:garding}). The sharpness of these estimates is measured by comparison against the short time asymptotics of Theorem \ref{barb:thm:30}.
\begin{theorem}\label{barb:thm:gb+ebd}{\rm \cite[Theorem 4.5]{barbatisdavies}}\par
Let $H$ be an operator of order $2m>n$ with real-valued coefficients in $L^{\infty}(\R^n)$. Assume that the principal coefficients
$\{a_{\ab}(x)\}_{|\alpha|=|\beta|=m}$ satisfy
\[
  \sum_{\subab}a^0_{\ab}v_{\alpha}\overline{v_{\beta}} \leq
  \sum_{\subab}a_{\ab}(x)v_{\alpha}\overline{v_{\beta}}
 \leq \mu \sum_{\subab}a^0_{\ab}v_{\alpha}\overline{v_{\beta}}, \quad v\in\C^{\nu}, \; x\in\R^n,
\]
for some $\mu\geq 1$, where $\{ a^0_{\ab}\}$ is a coefficient matrix for $(-\Delta)^m$. Then for any $\epsilon>0$ there exists
$c_{\epsilon}$ such that the heat kernel of $H$ satisfies
\[
| K(t,x,y)| < c_{\epsilon}t^{-\frac{n}{2m}}\exp\left\{ -( \rho(m,\mu)-\epsilon)\frac{|x-y|^{\frac{2m}{2m-1}}}{t^{\frac{1}{2m-1}}} + c_{\epsilon}t   \right\},
\]
for all $t\in\R_+$ and $x,y\in\R^n$, where
\[
\rho(m,\mu)=(2m-1)(2m)^{-2m/(2m-1)}\mu^{1/(2m-1)}\Big[ \sin\Big(\frac{\pi}{4m-2}\Big)^{-2m+1} +C\mu^m(\mu-1)\Big]^{-\frac{1}{2m-1}},
\]
and the constant $C$ depends only on $m$ and $n$. In particular $\rho(m,\mu)=\sigma_m +O(\mu-1)$ as $\mu\to 1^+$.
\end{theorem}
While Theorem \ref{barb:thm:gb+ebd} provides useful information when $H$ is close to $(-\Delta)^m$, it is clearly not very effective when $H$ is an arbitrary elliptic
operator. In such a case, the Finsler distance should play a role. Since definition (\ref{barb:def:fd}) is meaningless when $H$ has measurable
coefficients, an alternative definition is required, as was the case for second order operators.
Denoting by $A(x,\xi)$ the principal symbol of $H$ we define
\[
\cE =\{  \phi\in C^1(\R^n)  \, : \; A(x,\nabla\phi (x)) \leq 1  \mbox{ for almost all } x\in\R^n\}.
\]
For operators with smooth coefficients the Finsler distance $d(x,y)$ is then also given by
\begin{equation}\label{barb:eq:210}
d(x,y)=\sup\{ \phi(y)-\phi(x) \; : \; \phi\in\cE\} \; ;
\end{equation}
see \cite[Lemma 1.3]{agmon1}. Hence we use (\ref{barb:eq:210}) to define the Finsler distance when $H$ has measurable coefficients.
We note that a simple approximation argument shows that in the definition of $\cE$ we could have required that $\phi\in C^{\infty}(\R^n)$.
Given $M>0$ we also define
\[
\cE_M =\{  \phi\in C^m(\R^n)  \; : \; A(x,\nabla\phi (x)) \leq 1  \; , \;
| \nabla^k\phi(x)| \leq M \; , \mbox{ a.e. } x\in\R^n, \; 2\leq k\leq m\}
\]
and the Finsler-type distance
\neweq{fintype}
d_M(x,y)=\sup\{ \phi(y)-\phi(x)  : \; \phi\in\cE_M\}.
\endeq
So $d_{\infty}(x,y)=d(x,y)$, but for finite $M$ we have $d_M(x,y)\leq d(x,y)$ in general.

We finally define the following measure of regularity of the principal coefficients of $H$,
\[
q_A =\max_{\subab}\dist_{L^{\infty}(\R^n)}(a_{\ab}, W^{m-1,\infty}(\R^n)).
\]
In particular $q_A=0$ if the principal coefficients are uniformly continuous.

\begin{theorem}\label{thm1:barb}{\rm \cite[Theorem 1]{barbatis}}\par
Let $2m>n$. Assume that the principal symbol $A(x,\xi)$ is strongly convex. For any $M>0$ and $\epsilon>0$ there exists a constant $\Gamma_{\epsilon,M}$ such that the
heat kernel of $H$ satisfies
\begin{equation}
| K(t,x,y) |< \Gamma_{\epsilon,M}t^{-\frac{n}{2m}}\exp\left\{ -(\sigma_m -Cq_A-\epsilon)
\frac{  d_M(x,y)^{\frac{2m}{2m-1}}}{   t^{\frac{1}{2m-1}}   } + \Gamma_{\epsilon, M}t   \right\},
\label{barb:eq:10}
\end{equation}
for all $t\in\R_+$ and $x,y\in\R^n$.
\end{theorem}

The constant $C$ in \eq{barb:eq:10} depends only on $m$, $n$ and the constants in G{\aa}rding's inequality (\ref{barb:garding}).
In relation to the last theorem we mention the following open problems:

\begin{problem}
\label{bb1}
{\em Is the term $Cq_A$ necessary in (\ref{barb:eq:10})? Under what assumptions can it be removed?}
\end{problem}

\begin{problem}
\label{bb2}
{\em Is it possible to replace $d_M(x,y)$ by $d(x,y)$ in (\ref{barb:eq:10})? Under what assumptions can?}
\end{problem}

\begin{problem}
{\em What is the role of strong convexity in the above theorems? What are the best possible results if we do not assume the strong convexity?}
\end{problem}

\begin{problem}
{\em For operators with regular coefficients obtain sharp heat kernel estimates when $2m\leq n$.}
\end{problem}

A partial answer to Problems \ref{bb1} and \ref{bb2} is provided in the next theorem under additional assumptions on the principal coefficients. Of course,
the questions remain as to what is the best possible result for measurable coefficients. The proof of the theorem is geometric and consists in showing that
$d_M/d \to 1$ as $M\to +\infty$, uniformly in $x$ and $y$.
\begin{theorem}\label{barb:thm:approx_fm}{\rm \cite[Corollary 3]{barbatis1999}}\par
Let $H$ be an elliptic operator of order $2m>n$ whose principal symbol $A(x,\xi)$ is strongly convex, is $C^{m+1}$ with respect to $x$ and satisfies
$|\nabla^k_xA(x,\xi)|\leq c|\xi|^{2m}$, $0\leq k\leq m+1$. Assume further that the map
\[
(x,\xi)\mapsto A(x,\xi)^{\frac{1}{2m}}
\]
defines a Finsler metric on $\R^n$ in the sense that \eq{defFin}-\eq{Fin2} are satisfied. Then the heat kernel of $H$ satisfies the estimate
\begin{equation}
| K(t,x,y) |< c_{\epsilon}t^{-\frac{n}{2m}}\exp\left\{ -(\sigma_m -\epsilon)
\frac{  d(x,y)^{\frac{2m}{2m-1}}}{ t^{\frac{1}{2m-1}} } + c_{\epsilon}t   \right\},
\label{barb:eq:20}
\end{equation}
for any $\epsilon>0$ and all $t\in\R_+$ and $x,y\in\R^n$.
\end{theorem}

We next consider singular operators with unbounded coefficients.
Let $a_{\ab}(x)=a_{\beta,\alpha}(x)$, $|\alpha|=|\beta|= m$, be real-valued functions in $L^{\infty}_{\rm loc}(\R^n)$. We fix $s>0$ and assume that the
weight $a(x)=1+|x|^s$ controls the size of the matrix $\{a_{\ab}\}$ in the sense that
\[
  c^{-1}a(x)|v|^2 \leq
  \sum_{\subab}a_{\ab}(x)v_{\alpha}\overline{v_{\beta}}
 \leq  c a(x)|v|^2 , \quad v\in\C^{\nu}, \; x\in\R^n.
\]
We consider the elliptic operator
\[
Hu=(-1)^m\sum_{\substack{|\alpha|=m \\ |\beta|=m}}D^{\alpha}\{a_{\ab}D^{\beta}u\}
\]
on $L^2(\R^n)$, defined by means of a quadratic form similarly to the uniformly elliptic case; see \cite{barbatis1998} for details.
For $M>0$ we then define the set
\[
\cE_M=\{ \phi\in C^m(\R^n) \; : \; A(x,\nabla\phi(x))\leq 1 , \; |\nabla^k\phi|\leq M a(x)^{-k/2m}, \mbox{ a.e. }x\in\R^n, \; 2\leq k\leq m\}
\]
and the Finsler-type distance \eq{fintype}.
The weight $a(x)$ induces the weighted $L^{\infty}$-norm $\|u\|_{L^{\infty}_a(\R^n)}= \sup_{\R^n}(|u|/a)$
and more generally the weighted Sobolev spaces
\[
W^{k,\infty}_a(\R^n)=\{ u\in W^{k,\infty}_{\rm loc}(\R^n) : \; |\nabla^ju(x)|\leq ca(x)^{(2m-j)/2m}, \;  \mbox{ a.e. }x\in\R^n, \; 0\leq j\leq m-1\}
\]
We set
\[
q_A =\max_{\subab}\dist_{L_a^{\infty}(\R^n)}(a_{\ab}, W^{m-1,\infty}_a(\R^n)).
\]
\begin{theorem}\label{thm50:barb}{\rm \cite[Section 2]{barbatis1998} and \cite[Theorem 2.2]{barbatis2004}}\par
Assume that $n$ is odd, that $0<s<2m-n$ and that the principal symbol of $H$ is strongly convex.
Then for any $M>0$ and $\epsilon>0$ there exists a constant $\Gamma_{\epsilon,M}$ such that the heat kernel of $H$ satisfies
\begin{equation}
| K(t,x,y) |< \Gamma_{\epsilon,M}t^{-s}\exp\left\{ -(\sigma_m -cq_A-\epsilon)
\frac{  d_M(x,y)^{\frac{2m}{2m-1}}}{   t^{\frac{1}{2m-1}}   } +  \Gamma_{\epsilon,M}t  \right\},
\label{barb:eq:60}
\end{equation}
for all $t\in\R_+$ and $x,y\in\R^n$.
\end{theorem}

\begin{problem}
\label{sing_op}
{\em Find out what happens when $n$ is even.}
\end{problem}

We end this section presenting a theorem of Dungey \cite{dungey}
for powers of operators. Let $(X,d)$ be a metric space and $\mu$ be a positive Borel measure on $X$.
Assume that $X$ is of uniform polynomial growth, that is there exists $c>0$ and $D,D^*\in\N$ such that the volume $V(x,r)$ of any ball $B(x,r)$ satisfies
\[
\darrn{ c^{-1}r^D \leq V(x,r)\leq cr^D,}{ \mbox{ if }r\leq 1,}{c^{-1} r^{D^*} \leq V(x,r)\leq cr^{D^*},}{ \mbox{ if }r\geq 1.}
\]
Accordingly let
\[
V(r)=\darr{r^D,}{r\leq 1,}{r^{D^*},}{r\geq 1.}
\]
\begin{theorem}\label{barb:thm:dungey}{\rm \cite[Theorem 1]{dungey}}\par
Let $H$ be a non-negative self-adjoint operator on $L^2(X,d\mu)$.
Assume that the semigroup $e^{-Ht}$ has an integral kernel $K(t,x,y)$ which is continuous in $(x,y)$ for all $t\in\R_+$ and satisfies the Gaussian estimate
\[
|K(t,x,y)| < c_{\epsilon}V(t)^{-\frac{1}{2}}\exp\left\{ -\big(\frac{1}{4}-\epsilon\big)\frac{d(x,y)^2}{t}\right\} ,
\]
for any $\epsilon>0$ and all $t\in\R_+$ and $x,y\in X$. Then for any integer $m\geq 2$ the semigroup generated by $-H^m$ has an integral kernel $K_m(t,x,y)$
which satisfies the Gaussian estimate
\[
|K_m(t,x,y)| < c_{\epsilon}V(t)^{-\frac{1}{2m}}\exp\left\{ -(\sigma_m-\epsilon)\frac{d(x,y)^{\frac{2m}{2m-1}}}{t^{\frac{1}{2m-1}}}\right\} ,
\]
for any $\epsilon>0$ and all $t\in\R_+$ and $x,y\in X$.
\end{theorem}

\subsection{More on the heat kernel of the biharmonic operator}\label{heatker}

In particular situations, much more can be said about the kernels relative to \eq{general}.
In this subsection we collect a number of properties related to the heat kernel of the polyharmonic operator $(-\Delta)^m$. All the information about the heat kernel of $(-\Delta)^m$ is contained in the functions $f_{m,n}$, since (cf. \eq{fmn})
\[
K(t,x,y)=\alpha_{m,n}t^{-n/2m}f_{m,n}\Big( \frac{|x-y|}{t^{1/2m}} \Big).
\]

We specialize to the case $m=2$ and we give some hints on how to obtain the corresponding results in the higher order case $m\ge3$.
For simplicity, we denote $f_n=f_{2,n}$.\par
When $m=2$, \eq{linear2} becomes the Cauchy problem
\neweq{heat}
\left\{\begin{array}{ll}
u_t+\Delta^2u=0\quad & \mbox{in }\R^n\times\R_+\,,\\
u(x,0)=u_0(x)\quad & \mbox{in }\R^n\, ,
\end{array}\right.
\endeq
whereas the kernels defined in \eq{polyker} read
\neweq{b}
f_{n}(\eta)=\eta ^{1-n}\int_{0}^{\infty }e^{-s^{4}}(\eta s)^{n/2}J_{(n-2)/2}(\eta s)\,ds\ .
\endeq
These kernels obey the following  recurrence formula, see \cite{fgg}:
\begin{equation}  \label{recursion}
f'_n(\eta)=-\eta\, f_{n+2}(\eta)\qquad\mbox{for all }n\ge1.
\end{equation}
Moreover, thanks to Evgrafov-Postnikov \cite{ep} (see also \cite[(1.10)]{li}), we know that the kernels have exponential decay at infinity. More precisely,
define the constants
$$
\sigma=\frac{3\sqrt[3]{2}}{16}\ ,\qquad K_n=\frac{1}{(2\pi)^{n/2}}\, \frac{1}{\sqrt{3}\,\cdot 2^{(n-3)/3}}\ ,
$$
then, in any space dimension $n\ge1$, we have
\begin{equation}\label{modulusn}
f_{n}(\eta)=\frac{K_n}{\alpha_{2,n}\eta^{n/3}}\, \left\{\cos\Big(\sqrt{3}\, \sigma\, \eta^{4/3}-\frac{n\pi}{6}\Big)+O(\eta^{-4/3})\right\}\,
e^{-\sigma\eta ^{4/3}}\qquad \mbox{as }\eta \to\infty\ .
\end{equation}
In \cite{aar} one can find the definition of the Gamma function and the power series expansion of the Bessel function:
$$\Gamma(y)=\int_0^\infty e^{-s}\, s^{y-1}\, ds\quad(y>0)\, ,\qquad
J_\nu(y)=\sum_{k=0}^\infty\, \frac{(-1)^k(y/2)^{2k+\nu}}{k!\, \Gamma(k+\nu+1)} \quad(\nu > -1)\ ,$$
as well as further properties of $\Gamma$ and $J_\nu$. This allows to obtain the representation of $f_n$ through power series:

\begin{theorem}\label{powerseries} {\rm \cite[Theorem 2.1]{fgg}}\par\noindent
For any integer $j\ge1$, we have
\begin{equation}  \label{even}
f_{2j}(\eta)\ =\ \sum_{k=0}^\infty\ (-1)^k\ \frac{\Gamma\left(\frac{k+j}{2}\right)}{2^{2k+j+1}\, k!\ (k+j-1)!}\ \eta^{2k}\ .
\end{equation}
For any nonnegative integer $j$, we have
\begin{equation}  \label{odd}
f_{2j+1}(\eta)\ =\ \frac{2^j}{\sqrt{8\pi}}\ \sum_{k=0}^\infty\ (-1)^k\ \frac{
(k+j)!\ \Gamma\left(\frac{2k+2j+1}{4}\right)}{k!\ (2k+2j)!}\ \eta^{2k}\ .
\end{equation}
In particular, $f_n(0)>0$ for all $n$ and
\begin{equation*}
f_1(\eta)=\frac{1}{\sqrt{8\pi}}\sum_{k=0}^\infty\ (-1)^k\ \frac{\Gamma\left(\frac{2k+1}{4}\right)}{(2k)!}\ \eta^{2k}\ ,\qquad f_2(\eta)=\frac{1}{4}
\sum_{k=0}^\infty\ (-1)^k\ \frac{\Gamma\left(\frac{k+1}{2}\right)}{[2^k\, k!]^2}\ \eta^{2k}\ .
\end{equation*}
\end{theorem}

Using the properties of the Bessel functions, the following third order ODE for the function $f_n$ was derived in \cite[Theorem 2.2]{fgg}
for any integer $n\geq 1$:
\begin{equation}\label{3ord}
f_n'''(\eta)+\frac{n-1}{\eta}f_n''(\eta)-\frac{n-1}{\eta^2}f_n'(\eta)-\frac{\eta}{4}f_n(\eta)=0
\end{equation}
or, equivalently,
\neweq{ode1}
\left(\Delta f_n \right)'(\eta)=\frac{\eta}{4}f_n(\eta)\ .
\endeq

According to \eq{modulusn} the kernel $f_n(\eta)$, and hence the biharmonic heat kernel, has infinitely many sign changes as $\eta\to \infty$,
see also previous work by Bernstein \cite{bernstein} when $n=1$. We refer to \cite{gg2} for further (minor) properties concerning the behavior of the
kernels at some special points.\par
We now rescale the kernel $f_n$ and define the function:
\neweq{expressed}
v_\infty(y)=2^{n/2}\,\alpha_n \,f_n(\sqrt2\,|y|)=
2^{(n+2)/4}\,\alpha_n \,|y|^{1-n/2}\int_0^\infty e^{-s^4}\,s^{n/2}\,J_{(n-2)/2}(\sqrt{2}\,|y|\,s)\,ds\qquad\forall y\in\mathbb{R}^n
\endeq
where $\alpha_n$ is given by
\begin{equation*}
\alpha _{n}^{-1}=\omega _{n}\int_{0}^{\infty }r^{n-1}f_{n}(r)\,dr=\integ f_n(|x|)\,dx\ ;
\end{equation*}
here $\omega _{n}$ denotes the surface measure of the $n$-dimensional unit ball (so that $\omega_1=2$). Note that $\int_{\R^n}v_\infty(y)dy=1$.
Although the functions $v_\infty$ and $f_n$ are strictly related we maintain the double notation since, in our setting,
they play quite different roles; the former is a stationary solution to \eq{eqv} below, the latter is the biharmonic heat kernel.
We aim to study the moments of the function $v_\infty $ defined in \eq{expressed}. The prototype monomial in $\R^n$ is given by
\neweq{Pm}
P_\ell(y)=y^{\ell}:= \prod_{i=1}^n\,y_i^{\ell_i}\qquad\mbox{for }\ell=(\ell_1,...,\ell_n)\in\N^n
\endeq
and its degree is $|\ell|=\sum_i\ell_i$. Then we define the $P_\ell$-moment of $v_\infty $ by
\neweq{MPm}
\M_{P_\ell}:=\integ P_\ell(y)\,v_\infty (y)\,dy
\endeq
and we have

\begin{theorem}\label{generalmoments} {\rm \cite[Theorem 2]{gazzola}}\par\noindent
For any $\ell=(\ell_1,...,\ell_n)\in\N^n$ the following facts hold:
\begin{enumerate}
\item $\M_{\Delta^2P_\ell}=-\,|\ell|\,\M_{P_\ell}$,
\item if $|\ell|\not\in4\N$ or if at least one of the $\ell_i$'s is odd, then $\M_{P_\ell}=0$,
\item if $|\ell|\in8\N$ and all the $\ell_i$'s are even, then $\M_{P_\ell}>0$,
\item if $|\ell|\in8\N+4$ and all the $\ell_i$'s are even, then $\M_{P_\ell}<0$.
\end{enumerate}
\end{theorem}

We have so far considered moments having polynomials of $y$ as weights;
we now consider powers of $|y|$ which are polynomials only for even integer powers. For any $b>-n$ we define
\neweq{MMbb}
\M_b:=\integ |y|^b\,v_\infty (y)\,dy\ .
\endeq
Note that for $b>-n$ the above integral is finite since $|y|^{b}\,v_\infty (y)\thicksim v_\infty (0)\,|y|^{b}$ as $y\to0$ and $v_\infty $ has exponential
decay at infinity according to \eq{modulusn} and \eq{expressed}. If $P_\ell(y)=|y|^\ell$ for some $\ell\in2\N$, then $\M_\ell$ coincides with $\M_{P_\ell}$ as defined
in \eq{MPm}. We are again interested in the sign of these moments. The following result holds:

\begin{theorem}\label{signsmoments} {\rm \cite[Theorem 4]{gazzola}}\par\noindent
Assume that $n\ge1$ and that $b>-n$. Then
\begin{eqnarray*}
\M_b>0 \, , & \qquad\mbox{for all }\; & b\in(-n,2)\bigcup\Big(\bigcup_{k=0}^\infty(8k+6,8k+10)\Big)\ ,\\
\M_b=0 \, , & \qquad\mbox{for all }\; & b\in 4\N+2\ ,\\
\M_b<0 \, , & \qquad\mbox{for all }\; & b\in\bigcup_{k=0}^\infty(8k+2,8k+6)\ .
\end{eqnarray*}
\end{theorem}

When $b\in(-n,0]$, Theorem \ref{signsmoments} was first proved in \cite[Proposition 3.2]{fgg}. Theorems \ref{generalmoments} and
\ref{signsmoments} give further information about the sign-changing properties of the kernels $f_n$ (recall \eq{expressed}), and they better
describe how these infinitely many sign changes occur. They also show that the sign of the moments of $f_n$ {\em do not} depend on $n$.\par\medskip
We conclude this section by explaining how the just described properties of the biharmonic heat kernels can possibly be extended to higher order polyharmonic
kernels. First of all, we recall that \cite[Theorem 4.1]{ep} (see also \cite[(1.10)]{li}) gives the following generalization
to \eq{modulusn} in any space dimension $n\ge1$:
\neweq{polyker2}
f_{m,n}(\eta)=\frac{K_{m,n}}{\eta^{n(m-1)/(2m-1)}}\, \left\{\cos\left(a_m\eta^{2m/(2m-1)}-b_{m,n}\right)+O(\eta^{-2m/(2m-1)})\right\}\,
e^{-\sigma_m\eta ^{2m/(2m-1)}}
\endeq
as $\eta\to\infty$ for some (explicit) positive constants $K_{m,n}$ and $b_{m,n}$ depending on $m$ and $n$, and some (explicit) positive constants $\sigma_m$
and $a_m$ depending only on $m$.\par
Next, we suggest the following

\begin{problem} {\em Determine a power series representation of the kind of Theorem \ref{powerseries} for the higher order kernels $f_{m,n}$ ($m\ge3$)
defined in \eq{polyker}.}\end{problem}

To this end, by arguing as in \cite{fgg} and using \cite[Section 4.62]{aar}, it may be useful to notice that \eq{recursion} still holds, independently of $m$.
Moreover, the following $(2m-1)$-order differential equation holds:
\begin{equation}\label{ode2}
\Big(\Delta^{m-1}f_{m,n}\Big)'(\eta)=\frac{(-1)^m}{2m}\, \eta\, f_{m,n}(\eta)\qquad\mbox{for all }n\ge1.
\end{equation}
It is straightforward that \eq{ode2} coincides with \eq{ode1} if $m=2$, whereas it reduces to $f'(\eta)=-\frac{1}{2}\eta f(\eta)$
whenever $m=1$ (recall that in the latter case, the kernel $f$ is {\em independent} of $n$).\par
With these two identities, one obtains results similar to Theorem \ref{signsmoments}. In particular, one has
\neweq{Cnm}
C_{m,n,\beta}:=\omega_n\int_0^\infty\eta^{n-1-\beta}f_{m,n}(\eta)\, d\eta>0\qquad\mbox{for all integers }n\ge1\mbox{ and all }\beta\in[0,n)
\endeq
where $\omega_n$ denotes the measure of the unit ball in $\R^n$. The proof of \eq{Cnm} can be obtained following the same lines as
\cite[Proposition 3.2]{fgg}, see \cite{gg2}.

\begin{problem} {\em Prove the full extension of Theorem \ref{signsmoments} to the case of general $m\ge2$. What are the signs of
$C_{m,n,\beta}$ for all $\beta\in(-\infty,n)$? How do they depend on $m$?}\end{problem}

\section{The Fokker-Planck equation}\label{fokker}

In some situations it is convenient to transform \eq{linear2} into a Fokker-Planck-type equation. Let
$$R(t):=(2mt+1)^{1/2m}$$
so that $R(t)^{2m-1}\,R'(t)\equiv1$. Also put
\neweq{change}
u(x,t):=R(t)^{-n}\,v\left(\frac{x}{R(t)},\log R(t)\right).
\endeq
Then take $\tau=\log R(t)$ and $y=x/R(t)$. Some lengthy but straightforward computations show that $v=v(y,\tau)$ solves
\neweq{eqv}
\left\{\begin{array}{ll}
v_\tau+\LL v =0\qquad & \mbox{in }\R^n\times\R_+\,,\\
v(y,0)=u_0(y)\qquad & \mbox{in }\R^n\,,
\end{array}\right.
\endeq
where
\neweq{LL}
\LL v:=(-\Delta)^mv-\nabla\cdot(y\,v)\ .
\endeq
We recall here some properties of the operator $\LL$ defined in \eq{LL}. The most relevant one is that, contrary to the second
order heat equation, the operator $\LL$ is not self-adjoint: we refer to \cite[Section 3]{egkp2} for some properties of the adjoint operator $\LL^*$.
Let $\sigma_m>0$ be as in \eq{sigma}-\eq{polyker2} and, for any $a\in[0,\sigma_m)$, consider the function
\neweq{rho}
{\rho_a }(x)=e^{a\,|x|^{2m/(2m-1)}}\,,\qquad x\in\mathbb{R}^n
\endeq
so that, in particular, ${\rho_a }\equiv1$ if $a=0$. For any such function ${\rho_a }$ consider the space $L^2_\rhoa(\R^n)$, the weighted
$L^2$-space endowed with the scalar product and norm
\neweq{Lrho}
(u,v)_{L^2_\rhoa(\R^n)}=\int_{\R^n}{\rho_a }(x)\,u(x)\,\bar{v}(x)\,dx\ ,\qquad\|u\|_{L^2_\rhoa(\R^n)}^2=(u,u)_{L^2_\rhoa(\R^n)}\ .
\endeq
Clearly, if $a=0$ we have $L^2_\rhoa(\R^n)=L^2(\R^n)$. Together with the space $L^2_\rhoa(\R^n)$, we consider the weighted Sobolev space
$H^{2m}_\rhoa(\R^n)$ endowed with the scalar product
$$\langle u,v\rangle_{H^{2m}_\rhoa(\R^n) }=\int_{\R^n}{\rho_a }(x)\sum_{|\alpha|\le2m} D^\alpha u(x)\,D^\alpha\bar{v}(x)\,dx \ .$$
By \cite[Proposition 2.1]{egkp2} we know that $\LL$ is a bounded linear operator from $H^{2m}_\rhoa(\R^n)$ onto $L^2_\rhoa(\R^n)$.\par
We now wish to characterize the spectrum of $\LL$. In particular, the kernel of $\LL$ is nontrivial; any function in the kernel is a
stationary solution to \eq{eqv}. As for \eq{expressed}, we rescale the kernels $f_{m,n}$ by setting
\neweq{vinfty}
v_\infty(y)=C_{m,n}\,f_{m,n}\Big((2m)^{1/2m}\,|y|\Big), \qquad\forall y\in\mathbb{R}^n
\endeq
where $C_{m,n}>0$ is a normalization constant chosen in such a way that $\int_{\R^n}v_\infty(y)\, dy=1$;
note that $v_\infty\in\Sc$, where $\Sc$ is the space of smooth fast decaying functions:
\neweq{Sc}
\Sc:=\{w\in C^\infty(\R^n)\,:\,|x|^a\,D^\alpha w(x)\to0\mbox{ as }|x|\to\infty\mbox{ for all }a\ge0\ ,\ \alpha\in\N^n\}\ .
\endeq
In fact, there exists a unique stationary solution to \eq{eqv} which belongs to $\Sc$:

\begin{theorem}\label{spectrum} {\rm \cite[Theorem 2.1]{egkp2}}\par\noindent
Up to a multiplication by a constant, there exists a unique nontrivial stationary solution to \eq{eqv} which belongs to $\Sc$. This solution $v$ is
radially symmetric and, if we further assume that $\int_{\R^n}v(y)\,dy=1$, it is explicitly given by $v_\infty$ in \eq{vinfty}.\par
Moreover, the spectrum of $\LL$ coincides with the set of nonnegative integers, $\sigma(\LL)=\N$. Each eigenvalue
$\lambda\in\sigma(\LL)$ has finite multiplicity and the corresponding eigenfunctions are given by
$$D^\alpha v_\infty \qquad\mbox{for }|\alpha|=\lambda\in\N\ .$$
The set of eigenfunctions is complete in $L^2_\rhoa(\R^n)$ for any $a\in[0,\sigma_m)$.
\end{theorem}

This fundamental (and elegant) result certainly deserves more investigation. Consider the (normalized) projection operator $P_\rhoa$ defined by
\neweq{P0}
P_\rhoa\,w:=\left(\int_{\R^n}{\rho_a}\,w\,v_\infty\,dx\right)\,\frac{v_\infty}{\|v_\infty\|_{L^2_\rhoa(\R^n)}^2}\qquad\mbox{for all }w\in L^2_\rhoa(\R^n)\ .
\endeq
We recall two problems suggested in \cite{gazzola}.

\begin{problem} {\em  Prove the generalized Poincar\'e-type inequality
$$
\|u-P_\rhoa u\|_{L^2_\rhoa(\R^n)}^2\le(u,\LL u)_{L^2_\rhoa(\R^n)} \qquad\mbox{for all }u\in H^{2m}_\rhoa(\R^n)\ .
$$
Although from Theorem \ref{spectrum} we know that the least nontrivial eigenvalue of $\LL$ is 1, since $\LL$ is not
self-adjoint the above inequality is by far nontrivial. In particular, prove (or disprove) the following:
$$(u,\LL u)_{L^2_\rhoa(\R^n)}=\int_{\R^n}{\rho_a }(x)\,\bar{u}(x)\,\LL u(x)\,dx\,\ge\|u\|_{L^2_\rhoa(\R^n)}^2  \qquad\mbox{for all }u\in[\ker\LL]^\perp\ .$$
}\end{problem}

\begin{problem} {\em Determine the convergence rate in $L^p$ (for $1\le p<\infty$) of the solution to \eq{eqv} towards its projection onto the kernel,
that is, onto the space spanned by $v_\infty$.}\end{problem}

\section{Asymptotic behavior of the solution}\label{asymptoticbeh}

In this section we shed some light on the long-time behavior of solutions to \eq{heat}.
The asymptotic behavior is better seen in the Fokker-Planck equation. When $m=2$, \eq{eqv} reads
\neweq{eqv2}
\left\{\begin{array}{ll}
v_\tau+\Delta^2v-\nabla\cdot(y\,v) v =0\qquad & \mbox{in }\R^n\times\R_+\,,\\
v(y,0)=u_0(y)\qquad & \mbox{in }\R^n\,.
\end{array}\right.
\endeq
We now study the moments of the solution $v$ to \eq{eqv2}. Let $\mathcal S$ be as in \eq{Sc}, let $u_0\in\mathcal S$ and consider the solution $v$ to
\eq{eqv2}. Let $P_\ell$ be as in \eq{Pm} and consider the (time-dependent) map
$$M_{P_\ell,u_0 }(\tau):=\integ P_\ell(y)\,v(y,\tau)\,dy=\integ y^{\ell}\,v(y,\tau)\,dy\ .$$
Let $v_\infty$ be as in \eq{expressed} and let $\M_{P_\ell}$ be as in \eq{MPm}. We have

\begin{theorem}\label{any} {\rm \cite[Theorem 3]{gazzola}}\par\noindent
Assume that $u_0\in\Sc$ is normalized in such a way that
\neweq{normalized}
\integ u_0(y)\, dy=\integ v_\infty(y)\, dy=1
\endeq
and let $v$ denote the solution to \eq{eqv2}. For any $\tau\ge 0$, the following facts hold:
\begin{eqnarray*}
  & (i) & \mbox{ $M'_{P_\ell,u_0 }(\tau)=-\,M_{\Delta^2P_\ell,u_0 }(\tau)-|\ell|\,M_{P_\ell,u_0 }(\tau)$ for all $\ell\in\N^n$,}\\
  & (ii) &  \mbox{ $M_{P_\ell,u_0 }(\tau)=e^{-|\ell|\,\tau}\integ P_\ell(x)\, u_0(x)\, dx$ for all $|\ell|\le3$,}\\
  & (iii) & \mbox{ $\lim_{\tau\to\infty}M_{P_\ell,u_0 }(\tau)=\M_{P_\ell}$ for all $\ell\in\N^n$.}
\end{eqnarray*}
\end{theorem}

By combining Theorems \ref{generalmoments} and \ref{any}, we infer

\begin{corollary}
Assume that $u_0\in\Sc$ is normalized in such a way that \eq{normalized} holds and let $v$ denote the solution to \eq{eqv2}. Then
$$\lim_{\tau\to\infty}M_{P_\ell,u_0 }(\tau)\ \left\{\begin{array}{ll}
=0\quad & \mbox{if $|\ell|\not\in4\N$ or if at least one of the $\ell_i$'s is odd,}\\
>0\quad & \mbox{if $|\ell|\in8\N$ and all the $\ell_i$'s are even,}\\
<0\quad & \mbox{if $|\ell|\in8\N+4$ and all the $\ell_i$'s are even.}
\end{array}\right.$$
\end{corollary}

In the particular case where $|\ell|=2k$ and $P_\ell(y)=|y|^{2k}$ we may give a simple characterization of the moments of a solution to \eq{eqv2}.
Consider a solution $v$ to \eq{eqv2} with initial data $u_0\in\Sc$. For all $b\ge0$ let $\M_b$ be as in \eq{MMbb} and put
$$M_{b,u_0 }(\tau):=\integ|y|^b\,v(y,\tau)\,dy\ .$$
We then have

\begin{theorem}\label{momentsv} {\rm \cite[Theorem 5]{gazzola}}\par\noindent
Assume that $u_0\in\Sc$ is normalized in such a way that \eq{normalized} holds and let $v$ denote the solution to \eq{eqv2}.
Then for any $k\in\N$, $k\ge2$, the above defined functions satisfy the following ODE
\neweq{equadiffM}
M'_{2k,u_0 }(\tau)+2k\,M_{2k,u_0 }(\tau)=-\,2k\,(2k-2)\,(2k+n-2)\,(2k+n-4)\,M_{2k-4,u_0 }(\tau)\ .
\endeq
Moreover, for any $k\in\N$, we have
\neweq{limM}
\lim_{\tau\to+\infty}M_{2k,u_0 }(\tau)=\M_{2k}
\endeq
and the following explicit representation
\neweq{explicit}
M_{2k,u_0 }(\tau)=\sum_{j=0}^k a_j^k\,e^{-2j\tau}\, ,
\endeq
where $a_0^k=\M_{2k}$ and
\begin{eqnarray*}
(i)& & a_k^k=M_{2k,u_0 }(0)+2k\,(k-1)\,(2k+n-2)\,(2k+n-4)\sum_{j=0}^{k-2}\frac{a_j^{k-2}}{k-j}\ ,\\
(ii)& & a_{k-1}^k=0\quad \mbox{if }k\ge1\ ,\\
(iii)& & a_j^k=-\,\frac{2k\,(k-1)\,(2k+n-2)\,(2k+n-4)}{k-j}\,a_j^{k-2}\quad \mbox{if }k\ge2\mbox{ and }j=0,...,k-2\ .
\end{eqnarray*}
In $(i)$ we use the convention that $\sum_{j=0}^{k-2}=0$ if $k\le1$. \end{theorem}

Formula \eq{explicit} shows, for instance, that
$$M_{0,u_0 }(\tau)\equiv\integ u_0(y)\,dy\ ,\qquad M_{2,u_0 }(\tau)=e^{-2\tau}\integ|y|^2\,u_0(y)\,dy\ ,$$
$$M_{4,u_0 }(\tau)=-\,2n\,(n+2)\integ u_0(y)\, dy+e^{-4\tau}\integ\Big[|y|^4+2n\,(n+2)\Big]u_0(y)\, dy\ ,$$
\begin{eqnarray*}
M_{6,u_0 }(\tau) & = & -\,6\,(n+4)\,(n+2)\, e^{-2\tau}\integ|y|^2\,u_0(y)\, dy\\
\ & \ & +e^{-6\tau}\left(\integ|y|^6\,u_0(y)\, dy+6\,(n+4)\,(n+2)\integ|y|^2\,u_0(y)\, dy\right)\ .
\end{eqnarray*}

If $b\not\in2\N$ (so that $|y|^b$ is not a polynomial) we may still define the map $M_{b,u_0}$ and, for all $b\in[4,\infty)$, we obtain
$$M'_{b,u_0 }(\tau)+b\,M_{b,u_0 }(\tau)=-\,b\,(b-2)\,(b+n-2)\,(b+n-4)\,M_{b-4,u_0 }(\tau)\ .$$

Note that Theorems \ref{any} and \ref{momentsv} also hold in a weaker form if $u_0\in L^1(\R^n)$ and $|y|^a\,u_0\in L^1(\R^n)$ for some $a\ge4$.
In this case, the statements hold true under the additional restriction that $|\ell|\le a$. In particular, we have the following

\begin{corollary}\label{L1}
Assume that $(1+|y|^4)\,u_0\in L^1(\R^n)$ and that \eq{normalized} holds. If $v$ denotes the solution to \eq{eqv2}, then
$$\lim_{\tau\to+\infty}\int_{\R^n}|y|^4\,v(y,\tau)\,dy=\M_4<0\ .$$
\end{corollary}

\section{Positivity preserving property}\label{ppp}

Contrary to the second order heat equation, no general positivity preserving property
(ppp in the sequel) holds for the Cauchy problem \eq{general}. By ppp, we mean here that positivity of the initial datum $u_{0}$ implies positivity (in space and time) for the solution $u=u(x,t)$ of \eq{general}; this is of course equivalent to the kernel $K(t,x,y)$ being non-negative.

Nevertheless,
by exploiting the properties of the kernels, some restricted and somehow hidden versions of ppp can be observed for the fourth order parabolic equation
\begin{equation}
\left\{
\begin{array}{ll}
u_{t}+\Delta ^{2}u=0\qquad & \mbox{in }\R^n\times\R_+\ ,\\
u\left( x,0\right) =u_{0}(x)\qquad & \mbox{in }\mathbb{R}^{n}\ ,
\end{array}
\right.  \label{linear}
\end{equation}
where $n\geq 1$ and $u_{0}\in C^{0}\cap L^{\infty }\left( \mathbb{R}^{n}\right)$. In this section we recall several weakened versions of ppp for the problem (\ref{linear}). We start however with a theorem about the general problem \eq{general} which provides quantative information on the positivity of the heat kernel near the diagonal $\{x=y\}$.

\begin{theorem}\label{barb:thm:lower}{\rm \cite[Theorem 6]{davies4}}\par\noindent
Let $H$ be an homogeneous elliptic operator of order $2m>n$ acting on $L^2(\R^n)$. There exists constants $c_1,c_2>0$ such that
\begin{equation}
K(t,x,y) \geq c_1 t^{-\frac{n}{2m}}
\label{barb:eq:lower}
\end{equation}
for all $t\in\R_+$ and $x,y\in\R^n$ such that $|x-y|^{2m}\leq c_2t$.
\end{theorem}

Theorem \ref{barb:thm:lower} states that the solution $u=u(x,t)$ to \eq{general} when $u_0(x)=\delta_{\{x=z\}}$ (the Dirac delta distribution
at some $z\in\R^n$) satisfies $u(x,t)>0$ whenever $|x-z|^{2m}\leq c_2t$. Therefore, one expects that if the mass of $u_0$ is ``concentrated'' in some small region of
$\R^n$ then ppp holds, at least in some part of that region. This can be made precise for the simplified problem \eq{linear} on which we focus our
attention for the rest of this section.

\begin{theorem}\label{positive} {\rm \cite[Theorem 1]{gg1}}\par\noindent
Assume that $0\not\equiv u_0\ge0$ is continuous and has compact support in $\mathbb{R}^n$.
Let $u=u(x,t)$ denote the corresponding bounded strong solution of \eq{linear}. Then,\par
$(i)$ for any compact set $K\subset\mathbb{R}^n$ there exists $T_K=T_K(u_0)>0$ such that $u(x,t)>0$ for all $x\in K$
and $t\ge T_K$;\par
$(ii)$ there exists $\tau=\tau(u_0)>0$ such that for all $t>\tau$ there exists $x_t\in\mathbb{R}^n$ such that $u(x_t,t)<0$.
\end{theorem}

The trivial example $u_0\equiv1$ shows that, at least for statement $(ii)$, the compact support assumption cannot be
dropped. By Theorem \ref{positive} we see that {\em negativity for (\ref{linear}) exists in general and goes to infinity.}
Fine results concerning the validity of the eventual positivity property in presence of a source, may be found in \cite{berchio}.\par
It appears instructive to combine Theorem \ref{positive}
with the following energy conservation laws obtained in \cite[Corollary 1]{gazzola}: let $u_0\in L^1(\R^n)$ and let $u$ be the
solution to \eq{linear}; then, for all $t>0$ we have
\neweq{u1}
\int_{\R^n}u(x,t)\,dx=\int_{\R^n}u_0(x)\,dx\ ,
\endeq
\neweq{u2}
\frac{d}{dt}\int_{\R^n}u(x,t)^2\,dx=-\,2\int_{\R^n}|\Delta u(x,t)|^2\,dx\ .
\endeq
Denote by $u^+=\max\{u,0\}$ and $u^-=-\min\{u,0\}$ the positive and negative parts of a function $u$, so that $u=u^+-u^-$.
Theorem \ref{positive} states that if $u_0\in C^0(\R^n)$ has compact support and $0\not\equiv u_0\ge0$ in $\R^n$, then $u^-(x,t)\not\equiv0$
for all $t>0$. Moreover, \eq{u1} states that the map
$$t\mapsto\int_{\R^n}u(x,t)\, dx\qquad(t\ge0)$$
is constant and equals a strictly positive number. Hence,
$$\int_{\R^n}u^-(x,t)\, dx>\int_{\R^n}u^-(x,0)\, dx=0\quad\mbox{for all }t>0\, ,$$
$$\int_{\R^n}u^+(x,t)\, dx>\int_{\R^n}u^+(x,0)\, dx=\int_{\R^n}u_0(x)\, dx\quad\mbox{for all }t>0\, ;$$
here we use redundant notations ($u^+(x,0)=u_0^+(x)=u_0(x)$ and $u^-(x,0)=u_0^-(x)=0$) in order to emphasize the {\it strict inequalities} between
the mass of the positive (respectively, negative) part of the solution $u=u(x,t)$ and the the mass of the positive (respectively, negative) part of
initial datum $u_0$.\par
On the other hand, \eq{u2} states that
$$t\mapsto\int_{\R^n}u(x,t)^2\, dx\qquad(\tau\ge0)$$
decreases and, in particular, that
$$\int_{\R^n}u^+(x,t)^2\, dx<\int_{\R^n}u_0(x)^2\, dx=\int_{\R^n}u_0^+(x)^2\, dx\qquad(t>0)\ .$$
Summarizing, the $L^2$-norm of the positive part of the solution $u$ is {\bf smaller} than the $L^2$-norm of the positive part of the initial datum $u_0$,
whereas the $L^1$-norm of the positive part of the solution $u$ is {\bf larger} than the $L^1$-norm of the positive part of the initial datum $u_0$.

\begin{problem} {\em Prove the counterpart of Theorem \ref{positive} for \eq{linear2} (for any $m\ge2$) when $0\not\equiv u_0\ge0$ is continuous
and has compact support in $\mathbb{R}^n$.}\end{problem}

Next, we consider initial data $u_{0}$ which are not compactly supported and which display a given decay behavior as $|x|\rightarrow \infty$.
We fix some arbitrary $\beta\ge0$ and consider the functional set
\begin{equation*}
\cC_\beta:=\{g\in C^0(\mathbb{R}^n;\mathbb{R}_+):\, g(0)>0\, ,\ g(x)=o(|x|^\beta)
\mbox{ as }|x|\to\infty\}\ .
\end{equation*}

In a suitable class of initial data, a positivity result for the linear Cauchy problem~(\ref{linear}) holds:

\begin{theorem}\label{heat2} {\rm \cite[Theorem 1.1]{fgg}}\par
Let $\beta\ge0$ and let $g\in\cC_{\beta }$. Let
\neweq{kindu0}
u_{0}(x)=\frac{1}{g(x)+|x|^{\beta }}.
\endeq
Let $u=u(x,t)$ be the corresponding solution of (\ref{linear}) and $K\subset\mathbb{R}^{n}$ be a compact set.

(i) If $\beta<n$, then there exists $\widetilde{C}_{n,\beta}>0$ such that
\begin{equation*}
\lim_{t\to+\infty}\ t^{\beta/4}u(x,t)=\widetilde{C}_{n,\beta}\, ,
\end{equation*}
uniformly with respect to $x\in K$.

(ii) If $\beta \geq n$ and $g(x)\equiv 1$, then there exists $\widetilde{D}_{n,\beta }>0$ such that
\begin{equation}
\begin{array}{ll}
\displaystyle\lim_{t\rightarrow +\infty }\ t^{n/4}(\log
t)^{-1}\,u(x,t)=\widetilde{D}_{n,n}\quad & \mbox{if }\beta =n \\
\displaystyle\lim_{t\rightarrow +\infty }\ t^{n/4}\,u(x,t)=\widetilde{D}_{n,\beta }\quad
& \mbox{if }\beta >n\,,
\end{array}
\label{gequal1}
\end{equation}
uniformly with respect to $x\in K$.
\end{theorem}

\begin{problem} {\em By using \eq{Cnm}, prove the counterpart of Theorem \ref{heat2} for \eq{linear2} when $u_0$ is as in \eq{kindu0}.} \end{problem}

The constants $\widetilde{C}_{n,\beta }$ and $\widetilde{D}_{n,\beta }$ in Theorem \ref{heat2} \textit{do not}
depend on $K$. What does depend on $K$ is the \textquotedblleft speed of convergence\textquotedblright , namely
how fast $t^{\beta /4}u(x,t)-\widetilde{C}_{n,\beta }$ converges to $0$ (and similarly for $\widetilde{D}_{n,\beta}$).
Let us also mention that if $\beta \geq n$, then for any $g\in\cC_{\beta }$ (not necessarily constant) one
still has that $\lim_{t\rightarrow +\infty }\ t^{\beta /4}u(x,t)=+\infty $ uniformly with respect to $x\in K$.

\begin{remark} {\em The quantitative positivity result of Theorem~\ref{heat2} provides strong enough
information to be applied also to semilinear problems, see \cite{fgg,gg}. At a first glance, this appears somehow unexpected,
since the techniques connected with the proof of Theorem~\ref{heat2} seem to be purely linear.}\end{remark}

Theorem \ref{heat2} does not clarify whether the eventual positivity for solutions of (\ref{linear}) is global or only local.
Theorem \ref{positive} suggests that negativity for the solution of (\ref{linear}) always exists and shifts to infinity, provided
$\beta$ is sufficiently large.

\begin{problem} {\em Prove Theorem \ref{positive} $(ii)$ for any $u_0$ as in \eq{kindu0} for $\beta$ large enough.}\end{problem}

On the other hand, if $u_{0}\equiv 1$ then the solution of (\ref{linear}) is
$u\left( x,t\right) \equiv 1.$ This trivial example shows that if $\beta =0$,
presumably one has \textit{global} eventual positivity for (\ref{linear}).
At least in the case $n=1$, this is also true if $\beta $ is positive but sufficiently small:

\begin{theorem}\label{globalpositivity} {\rm \cite[Proposition A.6]{fgg}}\par\noindent
We assume that $n=1$ and $u_0 (x)=|x|^{-\beta}$. For $\beta>0$ sufficiently small,
the corresponding solution of (\ref{linear}) given by
$$u(x,t)=\alpha_n\int_{\mathbb{R}^{n}}\frac{f_{n}(|z|)}{|x-t^{1/4}z|^\beta}\,dz$$
is positive in $\R\times\R_+$.
\end{theorem}

\begin{problem} {\em Prove Theorem \ref{globalpositivity} in any space dimension $n\ge1$.}\end{problem}

By combining \eq{u1} with Corollary \ref{L1} and with Theorem \ref{globalpositivity}, we obtain

\begin{corollary}\label{L2}
Assume that $u_0>0$ a.e.\ in $\R^n$.\par
$(i)$ If $(1+|x|^4)\,u_0\in L^1(\R^n)$, then the solution $u$ to \eq{heat} changes sign.\par
$(ii)$ If $n=1$, there exists $\beta_0>0$ such that if $\beta\in(0,\beta_0)$ and $u_0(x)=|x|^{-\beta}$, then the solution $u$ to \eq{heat}
is a.e.\ positive in $\R\times\R_+$.
\end{corollary}

Corollary \ref{L2} can be interpreted as follows. From Theorem \ref{positive} we know that solutions $u$ to \eq{heat} with compactly supported nonnegative
initial data $u_0$ display the eventual local positivity property, that is, $u(x,t)$ becomes eventually positive on any compact subset of $\R^n$
but it is always strictly negative somewhere in a neighborhood of $|x|=\infty$. This happens because the biharmonic heat kernels exhibit oscillations and,
outside the support of $u_0$, they ``push below zero'' the initial datum. The same happens if $u_0>0$ but $u_0$ is ``very close to zero'',
see statement $(i)$. On the other hand, if $u_0>0$ and $u_0$ is ``far away from zero'' then the kernels do not have enough negative strength to push
the solution below zero, see statement $(ii)$. The trivial case $u_0\equiv1$ (which is a stationary solution to \eq{heat}!) well explains this situation.

Finally, the following result shows that in general, we cannot expect neither global positivity nor uniform bounds for eventual positivity.

\begin{theorem}\label{heat22} {\rm \cite[Theorem 1.2]{fgg}}\par\noindent
Let $\beta \in \left( 0,n\right) .$ For any $T>1$ there exists $g\in\cC_{\beta }$ such that if
\begin{equation*}
u_{0}(x)=\frac{1}{g(x)+|x|^{\beta }}
\end{equation*}
then, the corresponding solution $u=u\left( x,t\right)$ of (\ref{linear})
satisfies $u\left( x_{T},T\right) <0$ for some $x_{T}\in \mathbb{R}^{n}.$\
\end{theorem}

\begin{problem} {\em Extend Theorems \ref{globalpositivity} and \ref{heat22} to \eq{linear2} for any $m\ge2$.}\end{problem}

\end{document}